\documentclass[reqno,a4paper,11pt]{amsart}
\usepackage{amssymb}
\usepackage{float}
\usepackage{tikz}
\usepackage[dvips, lmargin=4cm, rmargin=4cm, tmargin=4cm, bmargin=4cm]{geometry}
\usepackage[colorlinks=true, citecolor=red, linkcolor=blue, urlcolor=blue]{hyperref}
\usepackage{mathrsfs}
\usepackage{cite}
\usepackage{xfrac}
\usepackage{graphicx}
\usepackage{hyperref}
\usepackage{amsmath}
\usepackage[T1]{fontenc}
\usepackage{enumitem}
\usepackage{dsfont}
\everymath{\displaystyle}
\newtheorem{theorem}{Theorem}[section]
\newtheorem{definition}{Definition}[section]
\newtheorem{remark}{Remark}[section]
\newtheorem{lemma}{Lemma}[section]
\newtheorem{proposition}{Proposition}[section]

\numberwithin{equation}{section}

\begin{document}
\usetikzlibrary{positioning}
\title[$(\Phi_x,\psi)$-Fractional Musielak-Orlicz Framework]{Fundamental Properties and Embedding Results in a Novel $\left(\Phi_x, \psi\right)$-Fractional Musielak Space with an Application to Nonlocal BVP
}
\author[Ayoub KASMI , El Houssine AZROUL , and Mohammed SHIMI ]
{Ayoub KASMI  $^1$,   El Houssine AZROUL $^2$,   and Mohammed SHIMI   $^3$ }
\address{Ayoub KASMI, Elhoussine AZROUL \newline
 Sidi Mohamed Ben Abdellah
 University,
 Faculty of Sciences Dhar el Mahraz, Laboratory of Modeling, Applied Mathematics, and Intelligent Systems, Fez, Morocco}
 \email{$^1$ayoub.kasmi@usmba.ac.ma}
 \email{$^2$elhoussine.azroul@usmba.ac.ma}
\address{Mohammed SHIMI\newline
 Sidi Mohammed Ben Abdellah University, ENS of Fez, Department of Mathematics and Computer Science, Laboratory of Mathematics and Applications to Engineering Sciences, Fez, Morocco.}
\email{$^3$mohammed.shimi2@usmba.ac.ma}
\subjclass[2010]{35R11, 46E30, 35S15, 35A15.}
\keywords{New $(\Phi_x,\psi)$-Fractional Musielak Spaces, Mountain pass theorem, Musielak space.}
\date{Month, Day, Year}
\begin{abstract}
In this paper, we introduce and study a novel class of generalized 
$(\Phi_x,\psi)$-fractional Musielak spaces 
$\mathcal{K}_{\Phi_x}^{\alpha, \beta, \psi}$, which extends classical fractional spaces 
and offers the flexibility to model heterogeneous and nonlinear phenomena with memory 
and nonlocal effects. A detailed and rigorous analysis of their functional structure is 
carried out. Several new properties and embedding results are established, highlighting 
the originality of the proposed framework and its relevance to nonlocal BVPs.\\
To illustrate the significance of this functional setting, we prove the existence of 
nontrivial solutions to a nonlinear fractional differential problem under an 
Ambrosetti--Rabinowitz type condition, using the mountain pass theorem.\\
Our results provide new perspectives for the analysis of nonlocal and nonhomogeneous 
equations in variable-exponent and Musielak--Orlicz settings.

\end{abstract}

\maketitle
\tableofcontents
\section{Introduction}\label{sec1}

Functional frameworks constitute one of the fundamental pillars of modern mathematical analysis, playing a central role in the rigorous formulation, theoretical investigation, and resolution of problems across various domains~\cite{Brezis,Adams}. These spaces not only provide a coherent mathematical structure, but also offer powerful tools for applications in optimization and PDE theory, with significant impact on various fields such as fluid mechanics and nonlinear elasticity signal processing, physics, and engineering~\cite{Treanţă,Romor,Růžička}. To truly understand complex deterministic or stochastic phenomena, however, it is insufficient to study only central tendencies such as the mean or median. One must also analyze variability, integrability, regularity, and long-term properties of solutions~\cite{Sobolev,Brezis}. For example, the Navier–Stokes equations for incompressible fluids require specific function spaces to ensure regularity, while variational problems in nonlinear elasticity involve spaces determined by the constitutive law of the material. \textit{This naturally raises the question of why new function spaces are developed when classical frameworks already exist.}

The answer lies in the ongoing generalization of functional spaces to meet increasingly complex and heterogeneous problems. Classical Lebesgue spaces $L^p(\Omega)$, measuring integrability via a fixed exponent $p$, underlie this hierarchy~\cite{Adams,Brezis}. Sobolev spaces $W^{k,p}(\Omega)$, which control weak derivatives up to order $k$, extend this, enabling the analysis of phenomena involving rates of change and leading to fundamental Sobolev embedding theorems~\cite{Adams}. Despite their success, Lebesgue and Sobolev spaces are inherently limited when modeling nonlinear or non-uniform phenomena, as they were originally designed for systems with uniform properties. Real-world applications, such as composite materials with position-dependent properties~\cite{Musielak,Diening}, require models in which the growth exponent $p(\cdot)$ varies with position, in contrast to traditional $L^p$ spaces. Variable exponents are particularly important because they can capture complex behaviors found in systems such as electrorheological fluids under varying electric fields~\cite{Abu-Jdayil1995,Abu-Jdayil1997}, vector fields in magnetostatics~\cite{Cekic2012}, or image restoration problems~\cite{Chen2006}. They also play a role in the mathematical modeling of quasi-Newtonian fluids and other nonlinear physical phenomena~\cite{Zhikov1997}. \\

Orlicz spaces $L^{\Phi}(\Omega)$, in which the classical power-law growth $|t|^p$ is replaced by a more general convex function $\Phi(\cdot)$, were introduced to extend the flexibility of Lebesgue spaces~\cite{Rao,Musielak,Mihăilescu}. These spaces allow for the treatment of problems where the growth is not necessarily polynomial. However, classical Orlicz spaces remain homogeneous in the sense that the growth function $\Phi(\cdot)$ is uniform throughout the domain $\Omega$. To overcome this limitation, Musielak-Orlicz spaces $L_{\Phi_x}(\Omega)$ allow the growth function $\Phi_x(\cdot)$ to vary with the position $x \in \Omega$~\cite{Musielak,Diening}, thereby enabling the modeling of systems whose nonlinearity may depend locally on the position. These spaces form a natural hierarchy classical Lebesgue spaces $L^p(\Omega)$ are contained in Orlicz spaces, which in turn are contained in Musielak-Orlicz spaces \cite{Musielak,Diening,Rao} namely
$$
L^p(\Omega) \subset L^{\Phi}(\Omega) \subset L_{\Phi_x}(\Omega).
$$
This hierarchy illustrates the progressive increase in flexibility and adaptability of functional frameworks.
\\

Fractional calculus has profoundly reshaped the modeling of nonlocal interactions in physical systems. In contrast to classical derivatives, which are local and describe the instantaneous rate of change at a point, fractional derivatives are integral operators whose evaluation inherently depends on the function over an extended domain. This built-in nonlocality enables the modeling of processes in which the present state is influenced by the entire past evolution of the system.

From a mathematical perspective, the kernels of fractional derivative operators serve as weighting functions that systematically account for past contributions. Such a property makes fractional calculus particularly effective for capturing hereditary behaviors and spatially distributed effects that are often neglected in classical models.

Overall, the adoption of fractional derivatives—including Riemann--Liouville, Caputo, and Hilfer formulations—provides a rigorous and unified framework for incorporating nonlocality and historical dependence into complex systems. This approach continues to expand the theoretical and applied analysis of nonlinear, heterogeneous, and memory-influenced media~\cite{Kilbas,Mainardi,Hilfer}.

Bringing these threads together, one may ask:\begin{center}
\texttt{what new mathematical structures emerge when variable growth conditions and fractional derivatives coexist within the same framework?} 
\end{center}
Combining Musielak-Orlicz spaces with fractional operators enables both spatial heterogeneity and memory effects to be modeled, allowing for the analysis of advanced materials and nonlocal equations.\\
The present work aims to introduce and develop the theory of generalized $(\Phi_x,\psi)$-Fractional Musielak Spaces $\mathcal{K}_{\Phi_x}^{\alpha, \beta, \psi}$. This construction unifies and significantly extends previous frameworks~\cite{kasmi,kasmi1,SRIVASTAVA,Jiao1}. The advantages of this approach are numerous:
\begin{itemize}
    \item[\checkmark]  Flexibility to model spatially varying nonlinearities
       \item[\checkmark] Applicability to a broad class of nonlinear and nonlocal equations
       \item[\checkmark] Compatibility with advanced analytical tools and techniques
       \item[\checkmark] Unification and extension of several previous functional frameworks
       \item[\checkmark] Relevance for fractional boundary value problems with variable growth conditions
\end{itemize}
To achieve this advantage, we had to overcome several difficulties and challenges, such as:

\begin{itemize}

\item \textit{Definition and Construction of the Spaces:}
the establishment of generalized $(\Phi_x,\psi)$-fractional Musielak spaces is particularly delicate, combining the Musielak--Orlicz modular structure with $\psi$-Hilfer fractional derivatives. The non-standard behavior of the $\Phi_x$ function led to substantial difficulties in obtaining suitable estimates and controlling the norms, which hinders the establishment of equivalences and essential functional properties for variational analysis.

\item \textit{ Nonlocality:} inherent to the generalized fractional operators $(\Phi_x,\psi)$, it introduces long-range interactions and memory effects. This makes classical local analytical techniques insufficient and necessitates substantial methodological adaptations.

\item \textit{Nonlinear Complexity:}
the strong nonlinearity of the operators requires advanced variational and topological techniques. Studying the existence, multiplicity, and stability of solutions remains a significant challenge.

\end{itemize}
To clarify the position and relevance of the new space $\mathcal{K}_{\Phi_x}^{\alpha, \beta, \psi}$, we provide below a schematic diagram illustrating the relationships between classical Lebesgue and Sobolev spaces, Orlicz and Musielak-Orlicz spaces, their fractional analogues, and the newly defined space. 
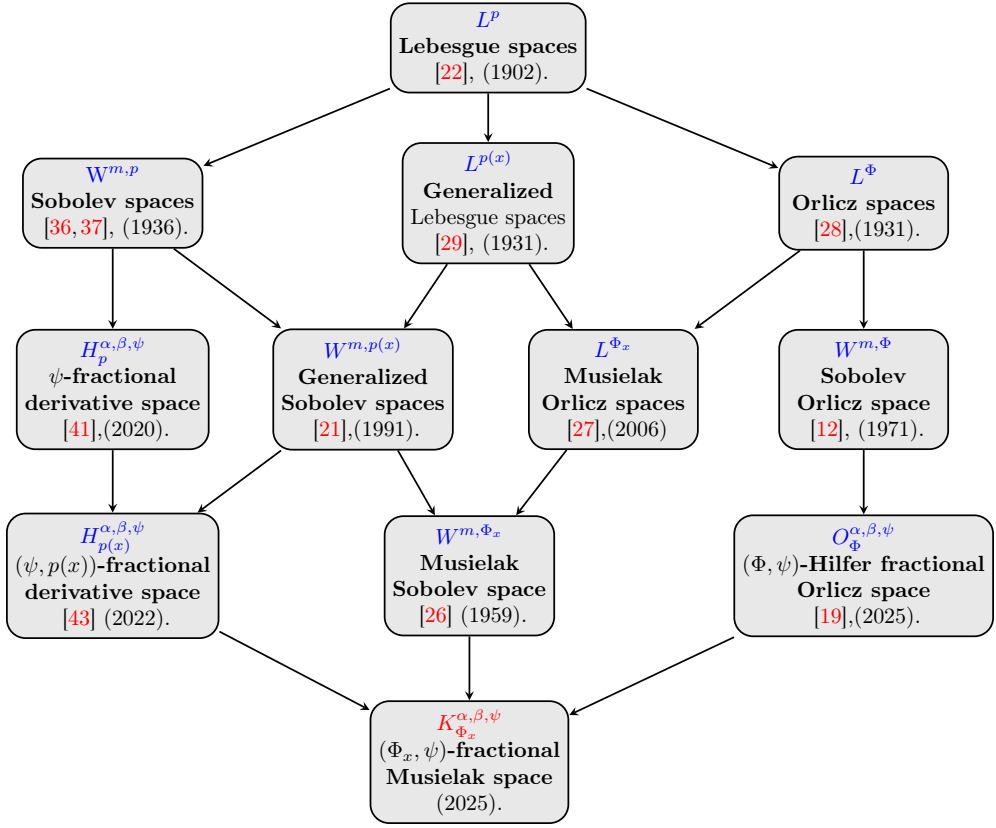
\begin{figure}[H]
\centering
\resizebox{1\textwidth}{!}{%
\begin{tikzpicture}[
  roundbox/.style={
    draw, rectangle, rounded corners=3mm,
    fill=gray!18,
    minimum width=2.7cm, minimum height=1cm,
    font=\small, align=center, thick
  },
  arrow/.style={thick,->,>=stealth}
]
\begin{scope}[shift={(-0.3,0)}]

\node[roundbox] (Lp) at (0,7.5) {\textbf{\textcolor{blue}{$L^p$}} \\ \textbf{Lebesgue spaces}\\ \cite{Lebesgue}, (1902).};

\node[roundbox] (Wmp)  at (-6,5) {\textcolor{blue}{W$^{m,p}$} \\ \textbf{Sobolev spaces}\\ \cite{Sobolev1,Sobolev}, (1936).};
\node[roundbox] (Lpx)  at (0,5)  {\textcolor{blue}{$L^{p(x)}$} \\\textbf{Generalized} \\Lebesgue spaces\\\cite{Orlicz0}, (1931).};
\node[roundbox] (Lphi) at (6,5)  {\textcolor{blue}{$L^{\Phi}$}\\ \textbf{Orlicz spaces}\\\cite{Orl},(1931).};

\node[roundbox] (Wmpx)  at (-2,2)  {$\textcolor{blue}{W^{m,p(x)}}$\\\textbf{Generalized}\\ \textbf{Sobolev spaces}\\
\cite{Kovacik},(1991).};
\node[roundbox] (Wphix) at (2,2)   {$\textcolor{blue}{L^{\Phi_x}}$\\\textbf{Musielak}\\\textbf{Orlicz spaces}\\ \cite{Musielak},(2006)};
\node[roundbox] (Wmphi) at (6,2)   {$\textcolor{blue}{W^{m,\Phi}}$\\\textbf{Sobolev}\\ \textbf{Orlicz space}\\\cite{Donaldson}, (1971). };

\node[roundbox] (Hmp)   at (-6,2) {$\textcolor{blue}{H_p^{\alpha, \beta,\psi}}$\\$\psi$\textbf{-fractional}\\ \textbf{derivative space}\\ \cite{A variational},(2020).};
\node[roundbox] (Hpx)   at (-6,-1) {$\textcolor{blue}{H^{\alpha, \beta, \psi}_{p(x)}}$\\ $(\psi,p(x))$\textbf{-fractional} \\\textbf{derivative space}\\ \cite{SRIVASTAVA} (2022).};
\node[roundbox] (OQGphi) at (6,-1) {$\textcolor{blue}{O^{\alpha, \beta, \psi}_{\Phi}}$\\$\left(\Phi,\psi\right)$\textbf{-Hilfer fractional}\\ \textbf{Orlicz space}\\ \cite{kasmi1},(2025).};
\node[roundbox] (Wmphi2) at (-0.3,-1) {$\textcolor{blue}{W^{m,\Phi_x}}$\\\textbf{Musielak}\\ \textbf{Sobolev space}\\ \cite{Musielak1} (1959).};
\node[roundbox] (KFinal) at (-0.3,-4) {$\textcolor{red}{K^{\alpha, \beta, \psi}_{\Phi_x}}$\\$\left(\Phi_x,\psi\right)$\textbf{-fractional}\\ \textbf{Musielak space}\\ (2025).};

% Flèches
\draw[arrow] (Lp) -- (Wmp);
\draw[arrow] (Lp) -- (Lpx);
\draw[arrow] (Lp) -- (Lphi);

\draw[arrow] (Wmp) -- (Wmpx);
\draw[arrow] (Wmp) -- (Hmp);
\draw[arrow] (Hmp) -- (Hpx);

\draw[arrow] (Lpx) -- (Wmpx);
\draw[arrow] (Lpx) -- (Wphix);
\draw[arrow] (Wmpx) -- (Hpx);
\draw[arrow] (Wmpx) -- (Wmphi2);

\draw[arrow] (Lphi) -- (Wphix);
\draw[arrow] (Lphi) -- (Wmphi); 
\draw[arrow] (Wphix) -- (Wmphi2);
\draw[arrow] (Wmphi) -- (OQGphi);

\draw[arrow] (Hpx) -- (KFinal);
\draw[arrow] (Wmphi2) -- (KFinal);
\draw[arrow] (OQGphi) -- (KFinal);

\end{scope}
\end{tikzpicture}%
} % fin de \resizebox

\vspace{0.5em}

\caption{Hierarchical diagram of functional spaces leading to the generalized $(\Phi_x,\psi)$-Fractional Musielak Spaces $\mathcal{K}_{\Phi_x}^{\alpha,\beta,\psi}$.}
\label{fig:spaces-diagram}
\end{figure}

This representation shows how each level of generalization brings additional flexibility, culminating in the broader framework proposed in this work.
\\

Our main results establish the analytical foundations of these spaces and highlight their importance in the study of fractional equations and nonhomogeneous models. In particular, we investigate a new class of non-integer (fractional) differential equations involving the $(\Phi_x, \psi)$-fractional operator, given by
\begin{equation} \label{P}
(\mathcal{P})~~
\left\{
\begin{array}{ll}
{}^{H}\mathcal{L}^{\alpha,\beta,\psi}_{\Phi_x} u
= h(x,u), & \text{in } \Lambda, \\[0.3cm]
u(0) = u(T) = 0, &
\end{array}
\right.
\end{equation}

 where : 
  \begin{itemize}[leftmargin=0.8cm]
 \item  $\alpha \in (0,1)$, $\beta \in [0,1]$ and $\Lambda=[0,T]$ .
 \item $a_x:\Lambda\times \mathbb{R}\to \mathbb{R}$ is  given in Section $\ref{sec2}.$
 \item ${}^{H}\mathcal{L}^{\alpha,\beta,\psi}_{\Phi_x} u={ }^{\mathbf{H}} \mathfrak{D}_{T}^{\alpha, \beta ; \psi}\!\left(
a_x\!\left({ }^{\mathbf{H}} \mathfrak{D}_{0^+}^{\alpha, \beta ; \psi} u\right)
{ }^{\mathbf{H}} \mathfrak{D}_{0^+}^{\alpha, \beta ; \psi} u
\right)$ is the fractional  $(\Phi_x,\psi)$-fractional operator.
 \item The nonlinearity $h:\Lambda \times \mathbb{R} \to \mathbb{R}$ is assumed to satisfy the following hypotheses:  
  \end{itemize}
\hypertarget{h1}{\textbf{(h$_1$)}} : The function $h$ is continuous on $\Lambda \times \mathbb{R}$, that is, $h \in \mathscr{C}(\Lambda \times \mathbb{R})$.\\
\hypertarget{h2}{\textbf{(h$_2$)}} : There exists a constant $\mu > k$, where $k$ denotes the constant introduced in the $\Delta_2$-condition in Section~\ref{sec2}, such that  
$$0 < \mu\, \mathcal{H}(t,u) \leq \, h(t,u)u, \quad \text{ for all } t \in \Lambda.$$

The assumption \hyperlink{h2}{\textbf{$($h$_2)$}} corresponds to the well-known Ambrosetti--Rabinowitz condition. 
This condition plays a fundamental role in the application of the Mountain Pass Theorem, 
a cornerstone result in critical point theory, which guarantees the existence of saddle-type critical points for nonlinear functionals.

In that context, Sousa et al. \cite{Sousa2022} established the existence of multiple solutions to the following fractional Dirichlet problem by utilizing the Nehari manifold approach:
$$
\small \left\{\begin{array}{l}
{ }^{\mathbf{H}} \mathbf{D}_T^{\alpha, \beta ; \psi}\left(\left|{ }^{\mathbf{H}} \mathbf{D}_{0+}^{\alpha, \beta ; \psi} u(x)\right|^{p(x)-2} { }^{\mathbf{H}} \mathbf{D}_{0+}^{\alpha, \beta ; \psi} u(x)\right)=\lambda a(x)|u|^{q(x)-2} u+b(x)|u|^{h(x)-2} u \quad \text { in } \Omega \\
u(x)=0 \quad \text { on } \partial \Omega.
\end{array}\right.
$$

For more details on this work, 
we refer the reader to~\cite{Sousa2022}. 
Notably, Lamine \textit{et al.}~\cite{Sousa2023} investigated the existence of solutions 
for a generalized fractional telegraph equation involving a class of $\psi$-Hilfer fractional derivatives 
combined with a $p(x)$-Laplacian type differential operator.
\begin{equation*}
\left\{
\begin{array}{l}
\varepsilon\, u_{tt} - { }^{\mathbf{H}} \mathbf{D}_T^{\alpha, \beta ; \psi}\left(\left|{ }^{\mathbf{H}} \mathbf{D}_{0^{+}}^{\alpha, \beta ; \psi} v\right|^{p(x)-2}{ }^{\mathbf{H}} \mathbf{D}_{0^{+}}^{\alpha, \beta ; \psi} v\right) + v_t = f(x, t), \\[0.3cm]
(x, t) \in \mathbf{Q}_T := \Omega \times (0, T).
\end{array}
\right.
\end{equation*}
\\

The structure of this paper is as follows. Section \ref{sec2} reviews the necessary theoretical background, including Musielak spaces, fractional calculus, and auxiliary results used throughout the work. In Section \ref{sec3}, we introduce the new $\left(\Phi_x,\psi\right)$-fractional Musielak space and establish several of its key properties. Section \ref{sec4} highlights the practical relevance of the proposed framework by applying the mountain pass theorem to analyze the fractional boundary value problem \hyperref[P]{($\mathcal{P}$)}.

\section{Basic Notation and Functional Background}\label{sec2}
In this section, we provide a comprehensive overview of the key results and foundational concepts in Musielak spaces and fractional calculus that will underpin the analyses and methods presented throughout this paper.

\subsection{ Musielak spaces }
This subsection is devoted to a concise presentation of the basic notions and main properties of Musielak spaces. More comprehensive discussions can be found in \cite{Musielak,Mihăilescu,Azroul}.
\begin{itemize}[leftmargin=0.5cm,parsep=0cm,itemsep=0cm,topsep=0cm]
\item Given $\Omega \subset \mathbb{R}^N$  an open set. 
Consider a function $a_x : \overline{\Omega} \times \mathbb{R} \to \mathbb{R}$ such that the mapping

$$
\varphi(x,t) = \varphi_x(t) :=
\begin{cases}
a(x,|t|)\, t, & \text{if } t \neq 0,\\[1mm]
0, & \text{if } t = 0,
\end{cases}
$$
is well-defined from $\overline{\Omega} \times \mathbb{R}$ into $\mathbb{R}$.

Moreover, we suppose that $\varphi$ fulfills the following property:  

\medskip
$(\varphi)$: For every $x \in \Omega$, the mapping $\varphi(x,\cdot): \mathbb{R} \to \mathbb{R}$ is an odd, strictly increasing homeomorphism from $\mathbb{R}$ onto $\mathbb{R}$.  
\end{itemize}
\medskip
\begin{definition}
Let $\Phi : \overline{\Omega} \times \mathbb{R} \to \mathbb{R}$ be defined by
$$
\Phi(x, t)=\Phi_x( t) = \int_0^t \varphi_x(s)\, ds, 
\qquad \forall x \in \overline{\Omega}, \ t \geq 0.
$$
We call $\Phi_x$ a \emph{Musielak function} provided that the following conditions are satisfied
\begin{itemize}
    \item[$(\Phi_1)$]: For every $x \in \Omega$, the mapping $\Phi(x,\cdot): [0,\infty) \to \mathbb{R}$ 
    is continuous, nondecreasing such that $\Phi(x,0)=0$ and $\Phi(x,t) > 0$ for all $t>0$, 
    with $\lim_{t\to \infty}\Phi(x,t)=\infty$.
    
    \item[$(\Phi_2)$]: For each $t \geq 0$, the mapping $\Phi(\cdot,t): \Omega \to \mathbb{R}$ is measurable.
\end{itemize}
\end{definition}

\begin{remark}
    The validity of condition $(\varphi)$ for $\varphi(x,\cdot)$ implies that $\Phi(x,\cdot)$ is a convex and nondecreasing mapping from $\mathbb{R}^+$ into $\mathbb{R}^+$.
\end{remark}

For the function $\Phi$ defined above, the Musielak class is defined as follows:
$$
K_{\Phi_x}(\Omega) = \Big\{ u : \Omega \to \mathbb{R} \ \text{measurable} \ : \ \int_{\Omega} \Phi_x(|u(x)|) \, dx < \infty \Big\},
$$
and the Musielak space (generalized Orlicz space),
$$
L_{\Phi_x}(\Omega) = \Big\{ u : \Omega \to \mathbb{R} \ \text{measurable} \ : \ \int_{\Omega} \Phi_x(\lambda |u(x)|) \, dx < \infty \ \text{for some } \lambda > 0 \Big\}.
$$

Equipped with the Luxemburg norm
\begin{equation}\label{NORM}
\|u\|_{\Phi_x} = \inf \Big\{ \lambda > 0 : \int_{\Omega} \Phi_x\left(\frac{|u(x)|}{\lambda}\right) \, dx \leq 1 \Big\},
\end{equation}
the space $L_{\Phi_x}(\Omega)$ becomes a Banach space.

The conjugate function of $\Phi_x$ is defined by
$$
\overline{\Phi}_x(t) = \int_0^t \overline{\varphi}_x(s) \, ds, \quad \forall x \in \overline{\Omega}, \ t \ge 0,
$$
where $\varphi_x : \mathbb{R} \to \mathbb{R}$ is given by
$$
\overline{\varphi}_x(t) = \sup \{ s \in \mathbb{R} : \varphi(x,s) \le t \}.
$$
Moreover, the following Hölder-type inequality holds (see \cite[Theorem 13.13]{Musielak}):
 
\begin{equation}\label{HOLDER}
    \left|\int_{\Omega} u v \mathrm{~d} x\right| \leqslant 2\|u\|_{\Phi_x}\|v\|_{\overline{\Phi}_x} \quad \text { for all } u \in L_{\Phi_x}(\Omega) \text { and } v \in L_{\overline{\Phi}_x}(\Omega) \text {. }
\end{equation}

In the sequel, we assume that
\begin{equation}\label{++}
    1<\varphi^{-}:=\inf _{t \geqslant 0} \frac{t \varphi_{x}(t)}{\Phi_{x}(t)} \leqslant \varphi^{+}:=\sup _{t \geqslant 0} \frac{t \varphi_{x}(t)}{\Phi_{x}(t)}<+\infty \quad \text { for all }x \in \overline{\Omega} \text{ and } t\geq 0.
\end{equation}
It follows from the above relation that $\Phi$ satisfies the $\Delta_2$-condition (see \cite{Mihăilescu}), that is 
\begin{equation}\label{delta}
    \Phi_x(2 t) \leq K \cdot \Phi_x(t), \quad \forall x \in \overline{\Omega}, t \geq 0,
\end{equation}
where $K$ is a positive constant.

Moreover, we suppose that $\Phi$ fulfills the following condition 
\begin{equation}\label{convex}
\text { for each } x \in \overline{\Omega} \text {, the function }[0, \infty) \ni t \rightarrow \Phi_x( \sqrt{t}) \text { is convex . }
\end{equation}
Relation \eqref{convex} guarantees that $L_{\Phi_x}(\Omega)$ is uniformly convex, and consequently reflexive (see \cite[Proposition 2]{Mihăilescu}).

The modular $\rho_{\Phi_x}$ associated with the Musielak space $L_{\Phi_x}(\Omega)$
$$ \rho_{\Phi_x}(u)=\int_{\Omega} \Phi_x(|u|) \mathrm{d} x$$
provides a natural measure of function size beyond the usual norm. Its strong link with the norm guarantees the equivalence of the induced topologies, allowing one to work interchangeably with the modular or the norm depending on the analytical context.

\begin{proposition}(\!\!\cite{Mihăilescu})\label{Pro21}
Assume that condition \eqref{++} is satisfied. Then, for every sequence $\left(u_n\right)$ and $u$ in $L_{\Phi_x}(\Omega)$, the following relations hold:
     \begin{enumerate}
         \item[i)]  $\|u\|_{\Phi_x}>1 \Rightarrow\|u\|_{\Phi_x}^{\varphi^-} \leqslant \rho_{\Phi_x}(u)\leqslant\|u\|_{\Phi_x}^{\varphi^{+}}$,
         \item[ii)]   $\|u\|_{\Phi_x}<1 \Rightarrow\|u\|_{\Phi_x}^{\varphi^+} \leqslant \rho_{\Phi_x}(u)\leqslant\|u\|_{\Phi_x}^{\varphi^{-}}$,
         \item[iii)] $\left\|u_n-u\right\|_{\Phi_x} \rightarrow 0 \quad \Leftrightarrow \quad \rho_{\Phi_x}(u_n-u) \rightarrow 0$.
     \end{enumerate}
\end{proposition}

%\begin{remark}
%If $\Phi(x,t) = t^p$ for some constant $p > 1$, then $L^{\Phi}(\Omega)$ coincides with the usual Lebesgue space $L^p(\Omega)$.  
%More generally, if $\Phi(x,t) = t^{p(x)}$, one recovers the variable exponent Lebesgue space $L^{p(x)}(\Omega)$.  
%\end{remark}

\subsection{$\psi$-Hilfer fractional derivative }  
	This subsection's remaining content is related to presenting the fractional Riemann-Lioville integral with respect to another function, the $\psi$-Hilfer
fractional derivative "\textbf{$\psi$-HFD}", and some results will be often used, (see \cite{Sousa}).
	\begin{itemize}
	    \item 
 Define $(a,b)$ as a non-empty interval in $\mathbb{R}$ (with $-\infty \leqslant a<b \leqslant \infty$). Consider $\psi$ a positive function defined on $(a,b)$ that is continuous and increasing, as well as $\psi \in \mathscr{C}^1(a,b)$. Given a real function $v$, we define the left and right-sided fractional integrals of a function $v$ with respect to $\psi$ on $[a, b]$ by
%  Consider a finite or infinite interval $-\infty \leq a<b \leq \infty$ on the real line $\mathbb{R}$, where $a<b$ and $n-1<\alpha<n$. Let $\psi$  a continuous, increasing, and positive function on $(a, b)$ with a continuous derivative $\psi^{\prime}$ on $(a, b)$. The left and right-sided fractional integrals of a function $v$ with respect to $\psi$ on $[a, b]$ are defined as follows:\\
		\begin{equation}\label{21}
			\mathbf{I}_{a+}^{\alpha ; \psi} v(x)=\frac{1}{\Gamma(\alpha)} \int_a^x \psi^{\prime}(t)(\psi(x)-\psi(t))^{\alpha-1} v(t) d t
		\end{equation}
		and
		\begin{equation}\label{22}
			\mathbf{I}_{b-}^{\alpha ; \psi} v(x)=\frac{1}{\Gamma(\alpha)} \int_x^b \psi^{\prime}(t)(\psi(t)-\psi(x))^{\alpha-1} v(t) d t.
		\end{equation}
	\\
	\item 
		Consider that $\psi^{\prime}(x) \neq 0$. The Riemann-Liouville derivatives of a function $v$ with respect to $\psi$ of order $\alpha \in (0,1)$, are defined by
		\begin{equation}\label{23}
			\begin{aligned}
				\mathfrak{D}_{a+}^{\alpha ; \psi} v(x) & =\left(\tfrac{1}{\psi^{\prime}(x)} \frac{d}{d x}\right)^n \mathbf{I}_{a+}^{n-\alpha ; \psi} v(x) \\
				& =\tfrac{1}{\Gamma(n-\alpha)}\left(\tfrac{1}{\psi^{\prime}(x)} \tfrac{d}{d x}\right)^n \int_a^x \psi^{\prime}(t)(\psi(x)-\psi(t))^{n-\alpha-1} v(t) d t
			\end{aligned}
		\end{equation}
		and
		\begin{equation}\label{24}
			\begin{aligned}
				\mathfrak{D}_{b-}^{\alpha ; \psi} v(x) & =\left(-\tfrac{1}{\psi^{\prime}(x)} \tfrac{d}{d x}\right)^n \mathbf{I}_{b-}^{n-\alpha ; \psi} v(x) \\
				& =\tfrac{1}{\Gamma(n-\alpha)}\left(-\tfrac{1}{\psi^{\prime}(x)} \tfrac{d}{d x}\right)^n \int_x^b \psi^{\prime}(t)(\psi(t)-\psi(x))^{n-\alpha-1} v(t) d t,
			\end{aligned}
		\end{equation}
		where $n=[\alpha]+1$.
	\\
	
		\item Under the aforementioned assumptions. The \textbf{$\psi$-HFD} left-sided and right-sided ${ }^{\mathbf{H}} \mathfrak{D}_{a+}^{\alpha, \beta ; \psi}(\cdot)$ and ${ }^{\mathbf{H}} \mathfrak{D}_{b-}^{\alpha, \beta ; \psi}(\cdot)$ of a function $v$ of order $\alpha \in (0,1)$ and type $\beta \in [0,1]$, are defined by
		\begin{equation}\label{25}
			{ }^{\mathbf{H}} \mathfrak{D}_{a+}^{\alpha, \beta ; \psi} v(x)=\mathfrak{D}_{a^+}=\mathbf{I}_{a+}^{\beta(n-\alpha) ; \psi}\left(\tfrac{1}{\psi^{\prime}(x)} \tfrac{d}{d x}\right)^n \mathbf{I}_{a+}^{(1-\beta)(n-\alpha) ; \psi} v(x)
		\end{equation}
		and
		\begin{equation}\label{26}
			{ }^{\mathbf{H}} \mathfrak{D}_{b-}^{\alpha, \beta ; \psi} v(x)=\mathfrak{D}_{b^-}=\mathbf{I}_{b-}^{\beta(n-\alpha) ; \psi}\left(-\tfrac{1}{\psi^{\prime}(x)} \tfrac{d}{d x}\right)^n \mathbf{I}_{b-}^{(1-\beta)(n-\alpha) ; \psi} v(x) .
		\end{equation}
		$(\ref{25})$ and $(\ref{26})$ can be expressed in the following way: 
		\begin{equation}\label{27}
			{ }^{\mathbf{H}} \mathfrak{D}_{a+}^{\alpha, \beta ; \psi} v(x)=\mathbf{I}_{a+}^{\eta-\alpha ; \psi} \mathfrak{D}_{a+}^{\eta ; \psi} v(x)
		\end{equation}
		and
		\begin{equation}\label{28}
			{ }^{\mathbf{H}} \mathfrak{D}_{b-}^{\alpha, \beta ; \psi} v(x)=\mathbf{I}_{b-}^{\eta-\alpha ; \psi} \mathfrak{D}_{b-}^{\eta ; \psi} v(x),
		\end{equation}
		with $\eta=\alpha(1-\beta)+n\beta$ and $\mathbf{I}_{a+}^{\eta-\alpha ; \psi}(\cdot),$  $\mathfrak{D}_{a+}^{\eta ; \psi}(\cdot),$ $\mathbf{I}_{b-}^{\eta-\alpha ; \psi}(\cdot),$ $\mathfrak{D}_{b-}^{\eta ; \psi}(\cdot)$ are defined in \eqref{21}, \eqref{22}, \eqref{23}, and \eqref{24}.
\\
    
	\item 
		If $v \in \mathscr{C}^n([a, b]),$ then
		$$
		\mathbf{I}_{a+}^{\alpha ; \psi}{ }^{\mathbf{H}} \mathfrak{D}_{a+}^{\alpha, \beta ; \psi} v(x)=v(x)-\sum_{k=1}^n \tfrac{(\psi(x)-\psi(a))^{\eta-k}}{\Gamma(\eta-k+1)} v_\psi^{[n-k]} \mathbf{I}_{a+}^{(1-\beta)(n-\alpha) ; \psi} v(a)
		$$
		and
		$$
		\mathbf{I}_{b-}^{\alpha ; \psi}{ }^{\mathbf{H}} \mathfrak{D}_{b-}^{\alpha, \beta ; \psi} v(x)=v(x)-\sum_{k=1}^n \tfrac{(-1)^k(\psi(b)-\psi(x))^{\eta-k}}{\Gamma(\eta-k+1)} v_\psi^{[n-k]} \mathbf{I}_{b-}^{(1-\beta)(n-\alpha) ; \psi} v(b).
		$$
\end{itemize}
We present here the main variational tool employed to establish our multiplicity result.
    \begin{theorem}\cite{Struwe} ({Mountain Pass Theorem })
    
Let $(X,\|\cdot\|)$ be a Banach space and let $\mathfrak{J}\in \mathscr{C}^1(X, \mathbb{R})$ satisfy the Palais-Smale (PS) condition. Suppose that $\mathfrak{J}(0)=0$ and that the following geometric hypotheses hold:
\begin{itemize}
    \item[$(G_1)$] There exist $L>0$ and $a>0$ such that $\mathfrak{J}(u) \geq a$ for all $u\in X$ with $\|u\|=L$;
    \item[$(G_2)$] There exists $u_0 \in X$ with $\|u_0\| > L$ such that $\mathfrak{J}(u_0) < 0$.
\end{itemize}
Then $\mathfrak{J}$ possesses a critical value $c \geq a$ which can be characterized by
\[
c := \inf_{\gamma\in\Gamma} \max_{t\in[0,1]} \mathfrak{J}(\gamma(t))
\]
where
\[
\Gamma := \left\{ \gamma \in C([0,1], X) \;\middle\vert\; \gamma(0)=0,\;\gamma(1)=u_0 \right\}.
\]
\end{theorem}

\section{New $\left(\Phi_x,\psi\right)$-fractional Musielak space}\label{sec3}
The new $\left(\Phi_x,\psi\right)$-fractional space, briefly presented in the introduction, is recalled and investigated in depth. A detailed and rigorous analysis of its structure is carried out, leading to the establishment of several new properties that highlight the originality of the proposed framework and its relevance to nonlocal problems.

Let $\alpha \in (0,1)$, $\beta \in [0,1]$, $\Lambda = [0,T]$ and let $\Phi_x$ be a Musielak function. We then introduce the left-sided $\left(\Phi_x,\psi\right)$-fractional Musielak space $\mathcal{K}^{\alpha,\beta,\psi}_{\Phi_x}(\Lambda, \mathbb{R})$ as follows
\begin{equation}\label{k}
        \mathcal{K}_{\Phi_x}:= 
\mathcal{K}^{\alpha,\beta,\psi}_{\Phi_x}=\left\{u\in L_{\Phi_x} (\Lambda, \mathbb{R}); \mathfrak{D}_{0^+}  u\in L_{\Phi_x}(\Lambda, \mathbb{R}) \right\}  .
        \end{equation}
This space is endowed with the norm given by
\begin{equation}\label{32}
    \|u\|_{\mathcal{K}_{\Phi_x}}=\|u\|_{\Phi_x}+[u]_{\Phi_x},
\end{equation}
where $[\cdot]_{\mathcal{K}_{\Phi_x}}$ being the Luxemburg norm defined as follows
$$
[u]_{\Phi_{x}}=\inf \left\{\delta>0: \int_{\Lambda} \Phi_{x}\left(\frac{\left|   \mathfrak{D}_{0^+}  u\right|}{\delta}\right) d x \leqslant 1\right\}.
$$
For any $v \in \mathcal{K}_{\Phi_x}^{\psi}$, we associate the modular 
$\rho_{\Phi_x(\cdot)}^{\psi} : \mathcal{K}_{\Phi_x}^{\psi} \to \mathbb{R}$ defined by
$$
\rho_{\Phi_{x}(\cdot)}^{\psi}(u) = \int_{\Lambda} \Big( \Phi_{x}(|u|) + \Phi_{x}(|\mathfrak{D}_{0^+} u|) \Big) \, dx.
$$
As a result, the norm $\|\cdot\|_{\mathcal{K}_{\Phi_x}^{\psi}}$ turns out to be equivalent to the modular norm
\begin{equation}\label{33}
 \|u\|_{\rho_{\Phi_x(\cdot)}^{\psi}}
 = \inf\left\{\delta > 0 : \rho_{\Phi_x(\cdot)}^{\psi}\!\left(\tfrac{u}{\delta}\right) \le 1\right\}.
\end{equation}
The space $\mathcal{K}^{0}_{\Phi_x}$ is defined as the completion of 
$\mathscr{C}_0^{\infty}(\Lambda,\mathbb{R})$ under the norm 
$\|\cdot\|_{\mathcal{K}_{\Phi_x}}$ introduced in~\eqref{32}, 
that is,
$$
\mathcal{K}^{0}_{\Phi_x}=\{\,u \in \mathcal{K}_{\Phi_x} : u(0)=u(T)=0\,\}.
$$

On the other hand, for any $v \in \mathcal{K}^{0}_{\Phi_x}$, we introduce the convex modular associated with the space $\mathcal{K}^{0}_{\Phi_x}$, which is defined by  
$$
{}^0\rho_{\Phi_x}^{\psi}(u)=\int_{\Lambda} \Phi_{x}(|\mathfrak{D}_{0^+} u|) \, dx .
$$

The norm associated with this modular, often referred to as the Luxemburg norm, is then given by  
$$
\|u\| = [u]_{\Phi_x}= \inf \left\{ \delta > 0 : \rho_{\Phi_{x}}^{\psi}\!\left(\frac{u}{\delta}\right) \leqslant 1 \right\}.
$$

\begin{remark}
We point out that the generalized fractional derivative space 
$\mathcal{K}_{\Phi_x}^{\alpha,\beta,\psi}$ introduced in~\eqref{k} 
represents a broad functional framework that unifies, under suitable choices 
of the function $a_x(\cdot)$ and the parameters $\alpha$, $\beta$, and $\psi$, 
several well-known fractional derivative spaces already studied in the literature. 
Indeed, by specifying these parameters appropriately, we recover a  wide range
of classical and recently developed spaces as particular cases, as illustrated below:

\begin{enumerate}
    \item If we take $a_x(t)=|t|^{p(x)-2}$ in $\eqref{k}$, we obtain the $(\psi,p(x))$-fractional derivative space
    $\mathcal{H}_{\kappa(x)}^{\gamma, \beta ; \chi}$, defined in~\textnormal{\cite{SRIVASTAVA}}, by
    $$
    \mathcal{H}_{p(x)}^{\gamma, \beta ; \chi}
    :=\left\{u \in L^{p(x)}(\Delta)\; ;\;
    {}^{\mathrm{H}}\mathcal{D}_{+}^{\gamma, \beta ; \chi}u \in L^{p(x)}(\Delta)
    \text{ and }u(\Delta)=0
    \right\}.
    $$
    
    \item If $\Phi_x(t)=\Phi(t)$ is independent of the variable $x$, we say that the space
    $\mathcal{K}_{\Phi_x}$ coincides with the space
    $\mathcal{O}_\Phi^{\alpha, \gamma, \psi}(\Lambda,\mathbb{R})$, defined in \cite{kasmi1}, by
    $$
    \mathcal{O}_\Phi^{\alpha, \gamma, \psi}
    :=\left\{ u \in L^{\Phi}(\Lambda,\mathbb{R})\; ;\;
    {}^{\mathbf{H}}\mathfrak{D}_{0^+}^{\alpha, \gamma ; \psi} u \in L^{\Phi}(\Lambda,\mathbb{R})
    \right\}.
    $$
    
    \item If $a_x(t)=|t|^{p-2}$ in $\eqref{k}$, the space $\mathcal{K}^{0}_{\Phi_x}$ reduce to the fractional derivative space
    $\mathbb{H}_p^{\alpha, \gamma ; \psi}$, defined in~\cite{A variational}, by
    $$
    \mathbb{H}_p^{\alpha, \gamma ; \psi}
    :=\left\{
    \begin{array}{l}
    u \in L^p(\Lambda,\mathbb{R}) \; ; \;
    {}^{\mathbf{H}}\mathbf{D}_{0+}^{\alpha, \gamma ; \psi} u \in L^p(\Lambda,\mathbb{R}),\\[0.3em]
    \mathbf{I}_{0+}^{\gamma(\gamma-1)} u(0)=0,\;
    \mathbf{I}_{T-}^{\gamma(\gamma-1)} u(T)=0
    \end{array}
    \right\}.
    $$
    \item If we take $a_x(t)=|t|^{p-2}$, $\psi(t)=t$, and $\beta \to 1$ in~\eqref{k},
    we recover the fractional derivative space $E_0^{\alpha, p}$, defined in \cite{Zhou}, by
    $$
    E_0^{\alpha, p}
    :=\left\{ u \in L^{p}([0,T]) \; ; \;
    {}_0^{c}D_t^{\alpha} u \in L^{p}([0,T]) \text{ and } u(0)=u(T)=0 \right\}.
    $$
\end{enumerate}

\end{remark}

%The diagram below provides a visual interpretation of the relationships among the various fractional derivative
%function spaces presented above. The generalized space 
%$\mathcal{K}_{\Phi_x}^{\alpha,\beta,\psi}(\Lambda)$, shown in red, serves as a unifying framework
%that encompasses several well-known spaces as particular cases. 
%In particular, when the function $a_x(t)$ and the parameters $\alpha$, $\beta$, and $\psi$ are specified appropriately,
%one can recover the spaces 
%$\mathcal{H}_{\kappa(x)}^{\gamma, \beta ; \chi}(\Lambda)$, 
%$\mathcal{O}_G^{\alpha, \gamma, \psi}(\Lambda)$, 
%and $\mathbb{H}_p^{\alpha, \gamma ; \psi}(\Lambda)$.
%This schematic representation highlights how the generalized Musielak–Orlicz framework
%naturally extends and connects these existing models within a single functional setting.

%\begin{tikzpicture}[scale=0.8]

% === Titre global ===
%\node[font=\bfseries] at (0,4.6) {Connections between Fractional Derivative Function Spaces}; % 

% Rectangle englobant
%\draw[thick] (-6,-3.5) rectangle (6,4);

% Cercles plus grands (rayon 3)
%\draw[blue,thick] (-1.9,0) circle (3);
%\draw[blue,thick] ( 1.9,0) circle (3);

%% Étiquettes des cercles
%\node at (-2.5,0) {$\mathcal{H}_{\kappa(x)}^{\gamma, \beta %; \chi}(\Lambda)$};
%\node at ( 2.5,0) {$\mathcal{O}_G^{\alpha, \gamma, \psi}(\Lambda)$};

% Intersection
%\node at (0,-0.2) {$\mathbb{H}_p^{\alpha, \gamma ; \psi}(\Lambda)$};

%% Étiquette du grand ensemble
%\node at (0,3.2) {$\textcolor{red}{\mathcal{K}_{\Phi_x}(\Lambda)}$};

%\end{tikzpicture}
\begin{theorem}
Let $\alpha \in (0,1)$, $\beta \in [0,1]$, and let $\Phi_x$ be a Musielak function defined on $\Lambda$.  
Then the space $(\mathcal{K}_{\Phi_x}, \|\cdot\|_{\mathcal{K}_{\Phi_x}})$ is a Banach space. 
Moreover, it is separable (respectively reflexive) if and only if $\Phi_x \in \Delta_2$ 
(respectively $\Phi_x, \overline{\Phi}_x \in \Delta_2$).  
In addition, if $\Phi_x \in \Delta_2$ and the function $t \mapsto \Phi_x(\sqrt{t})$ is convex, 
then $\mathcal{K}_{\Phi_x}$ is uniformly convex.
\end{theorem}

\noindent \textbf{\textit{Proof.}}
We define the operator $\mathfrak{T}$ by  
$$\mathfrak{T}: \mathcal{K}_{\Phi_x} \longrightarrow L_{\Phi_x}(\Lambda) \times L_{\Phi_x}(\Lambda) = P, \quad
u \longmapsto \left(u, \mathfrak{D}_{0^+} u\right).$$
Observe that for any $v \in \mathcal{K}_{\Phi_x}$,
$$\|\mathfrak{T}(u)\|_{P} = \|u\|_{\Phi_x} + \|\mathfrak{D}_{0^+} u\|_{\Phi_x} 
= \|u\|_{\Phi_x} + [u]_{\Phi_x} = \|u\|_{\mathcal{K}_{\Phi_x}}.$$
Hence, $\mathfrak{T}$ is an isometry. Because $P$ is a Banach space \cite{Adams}, it follows that 
$\mathcal{K}_{\Phi_x}$ is also a Banach space. Moreover, as $L_{\Phi_x}(\Lambda)$ is separable and uniformly convex 
(and therefore reflexive) (see \cite{Musielak,Mihăilescu}), we deduce that $\mathcal{K}_{\Phi_x}$ inherits these 
properties and is reflexive, separable, and uniformly convex.
\hfill$\Box$

%*****************************************************
\begin{proposition}
Assume that condition~\eqref{++} holds. Then, for all $x \in \Lambda$ and all $t \geq 0$, 
\[
\overline{\Phi}_{x}\big(\varphi_{x}(t)\big) \leq \varphi^{+} \, \Phi_{x}(t).
\]
\end{proposition}

\noindent \textbf{\textit{Proof.}}
Let us recall that
$$
\overline{\varphi}_{x}(t)
= \sup \{ s : \varphi_x(s)\le t \}, \qquad
\overline{\Phi}_{x}(t)
= \int_{0}^{t} \overline{\varphi}_{x}(s)\, ds,
\quad \forall x\in \Lambda,\ t\ge0.
$$

Furthermore, for every $x\in \Lambda$, the function
$\varphi_{x}:\mathbb{R}\to \mathbb{R}$ is an increasing homeomorphism,
in particular from $\mathbb{R}^{+}$ onto $\mathbb{R}^{+}$.
Hence, for each $x$, the map $t\mapsto \varphi_{x}(t)$ has an inverse
function $t\mapsto \varphi_{x}^{-1}(t)$.
It follows that $\varphi_{x}(s)\le t$ if and only if
$s\le \varphi_{x}^{-1}(t)$, which implies
$$
\overline{\varphi}_{x}(t) = \varphi_{x}^{-1}(t).
$$
Consequently,
$$
\overline{\Phi}_{x}(t)
= \int_{0}^{t} \varphi_{x}^{-1}(s)\, ds,
\quad \forall x\in \Lambda,\ t\ge0.
$$
Next, using a change of variables in the integral, we have
$$
\overline{\Phi}_{x}\bigl(\varphi_{x}^{-1}(s)\bigr)
= \int_{0}^{\varphi_{x}^{-1}(s)} \varphi_{x}(\theta)\, d\theta
= \int_{0}^{s} r \,\frac{d}{dr}\varphi_{x}^{-1}(r)\, dr
= s\,\varphi_{x}^{-1}(s) - \overline{\Phi}_{x}(s),
\quad \forall x\in \Lambda,\ s\ge0.
$$
This identity yields
$$
\overline{\Phi}_{x}(s) \le s\,\varphi_{x}^{-1}(s),
\quad \forall x\in \Lambda,\ s\ge0.
$$
Choosing $s=\varphi_{x}(t)$ gives
$$
\overline{\Phi}_{x}\bigl(\varphi_{x}(t)\bigr)
\le t\,\varphi_{x}(t),
\quad \forall x\in \Lambda,\ t\ge0.
$$
Finally, by assumption $\eqref{++}$,
$$
t\,\varphi_{x}(t) \le \varphi^{+}\,\Phi_{x}(t).
$$
By combining the above inequalities, we deduce that
\[
\overline{\Phi}_{x}\bigl(\varphi_{x}(t)\bigr)
\le \varphi^{+}\,\Phi_{x}(t),
\quad \text{for all } x \in \Lambda \text{ and } t \ge 0,
\]
which concludes the proof.
 \hfill$\Box$
%$********************************************
 \begin{proposition}
Let $\Phi_x$ be a Musielak function. 
Then for almost every $x$ in $\Lambda$ and for every $s\ge0$ one has
\begin{equation}\label{S}
\Phi_x(s)\le s\,\varphi_x(s)\le \Phi_x(2s).
\end{equation}
\end{proposition}
\noindent \textbf{\textit{Proof.}}
The case $s=0$ is immediate. Assume $s>0$.

Since $\varphi_x$ is nondecreasing, for every $t\in[0,s]$ we have $\varphi_x(t)\le\varphi_x(s)$, therefore
$$
\Phi_x(s)=\int_0^s\varphi_x(t)\,dt \le \int_0^s\varphi_x(s)\,dt = s\,\varphi_x(s),
$$
which proves the left inequality.

On the other hand, for every $t\in[s,2s]$ we have $\varphi_x(t)\ge\varphi_x(s)$, hence
$$
\Phi_x(2s)=\int_0^{2s}\varphi_x(t)\,dt \ge \int_s^{2s}\varphi_x(t)\,dt \ge \int_s^{2s}\varphi_x(s)\,dt = s\,\varphi_x(s).
$$
This yields the right inequality. Combining both estimates gives the claim. \hfill$\Box$\\
%**************************************************
\begin{proposition}\label{P32}
Let $v \in \mathcal{K}_{\Phi_x}$ and assume that condition~\eqref{++} holds. 
Then the following inequalities are satisfied:

\begin{equation}\label{34}
\text{If } [u]_{\Phi_x} > 1, \quad 
[u]_{\Phi_x}^{\varphi^{-}} \le {}^0 \rho_{\Phi_x}^{\psi}(u) \le [u]_{\Phi_x}^{\varphi^{+}} \quad \text{for all } u \in \mathcal{K}_{\Phi_x},
\end{equation}

\begin{equation}\label{35}
\text{If } [u]_{\Phi_x} < 1, \quad 
[u]_{\Phi_x}^{\varphi^{+}} \le {}^0 \rho_{\Phi_x}^{\psi}(u) \le [u]_{\Phi_x}^{\varphi^{-}} \quad \text{for all } u \in \mathcal{K}_{\Phi_x}.
\end{equation}
\end{proposition}

\noindent \textbf{\textit{Proof.}}
First, we establish that  
$$
{}^0\rho_{\Phi_x}^{\psi}(u) \leq [u]_{\Phi_x}^{\varphi^{+}}, 
\quad \text{for all } u \in \mathcal{K}_{\Phi_x} \text{ with } [u]_{\Phi_x} > 1.
$$  
Indeed, since  
$$
\varphi^{+} \geq \frac{s \varphi(s)}{\Phi_x(s)}, 
\quad \forall\, s>0,
$$  
it follows that, for any $r>1$, we have  
$$
\log \big(\Phi_x(r s)\big) - \log \big(\Phi_x(s)\big) 
= \int_s^{r s} \frac{\varphi(\tau)}{\Phi_x(\tau)}\, d\tau 
\leq \int_s^{r s} \frac{\varphi^{+}}{\tau}\, d\tau 
= \log \left(r^{\varphi^{+}}\right).
$$  
Hence, we deduce that  
\begin{equation}\label{36}
\Phi_x(r s) \leq r^{\varphi^{+}} \Phi_x(s), 
\quad \forall\, s>0, \; r>1.
\end{equation}  
Next, let $u \in \mathcal{K}_{\Phi_x}$ with $[u]_{\Phi_x}>1$. Using \eqref{36} together with the definition of the Luxemburg norm, we obtain  
$$
\begin{aligned}
\int_{\Lambda} \Phi_x\!\left(\left|\mathfrak{D}_{0^+} u\right|\right) dx 
&= \int_{\Lambda} \Phi_x\!\left([u]_{\Phi_x}\, \tfrac{|\mathfrak{D}_{0^+} u|}{[u]_{\Phi_x}}\right) dx \\
&\leq [u]_{\Phi_x}^{\varphi^{+}} \int_{\Lambda} \Phi_x\!\left(\tfrac{|\mathfrak{D}_{0^+} u|}{[u]_{\Phi_x}}\right) dx \\
&\leq [u]_{\Phi_x}^{\varphi^{+}}.
\end{aligned}
$$
Now, we prove that  
$$
[u]_{\Phi_x}^{\,\varphi^{-}} \leq {}^0\rho_{\Phi_x}^{\psi}(u), 
\quad \text{for all } u \in \mathcal{K}_{\Phi_x} \text{ with } [u]_{\Phi_x} > 1.
$$  
By employing a similar argument as in relation \eqref{36}, we deduce that  
\begin{equation}\label{38}
\Phi_x(r s) \geq r^{\varphi^{-}} \Phi_x(s), 
\quad \forall\, s>0 \text{ and } r>1 .
\end{equation}  
Let $u \in \mathcal{K}_{\Phi_x}$ with $[u]_{\Phi_x} > 1$.  
Choose $\sigma \in \left(1,[u]_{\Phi_x}\right)$.  
Since $\sigma < [u]_{\Phi_x}$, the definition of the Luxemburg norm yields  
$$
\int_{\Lambda} \Phi_x\!\left(\tfrac{\left|\mathfrak{D}_{0^+} u\right|}{\sigma}\right) dx > 1,
$$  
otherwise this would contradict the definition of the norm.  

Therefore, using \eqref{38}, we obtain  
$$
\begin{aligned}
\int_{\Lambda} \Phi_x\!\left(\left|\mathfrak{D}_{0^+} u\right|\right) dx
&= \int_{\Lambda} \Phi_x\!\left(\beta \cdot \tfrac{\left|\mathfrak{D}_{0^+} u\right|}{\sigma}\right) dx \\
&\geq \sigma^{\varphi^{-}} \int_{\Lambda} \Phi_x\!\left(\tfrac{\left|\mathfrak{D}_{0^+} u\right|}{\sigma}\right) dx \\
&\geq \sigma^{\varphi^{-}}.
\end{aligned}
$$  

Finally, for $\sigma \to [u]_{\Phi_x}$, We therefore deduce that inequality~\eqref{34} is satisfied.  

Next, we show that 
$$
{}^0 \rho_{\Phi_x}^{\psi}(u) \le [u]_{\Phi_x}^{\varphi^{-}}, 
\quad \text{for all } u \in \mathcal{K}_{\Phi_x} \text{ with } [u]_{\Phi_x} < 1.
$$
Using an argument analogous to that employed in~\eqref{36}, we obtain
\begin{equation}\label{39}
    \Phi_x(s) \leqslant \tau^{\varphi^{-}} \Phi_x\left(\frac{s}{\tau}\right) \text { for all } s>0 \text { and } \tau \in(0,1) \text {. }
\end{equation}
Let $u \in \mathcal{K}_{\Phi_x}$ with $[u]_{\Phi_x}<1$. From the definition of the norm define by $(\ref{NORM})$ and the relation $(\ref{39})$, we conclude  
$$
\begin{aligned}
 \int_{\Lambda} \Phi_x\left(\left|\mathfrak{D}_{0^+} u\right|\right) d x & =\int_{\Lambda} \Phi_x\left(\tfrac{\left|\mathfrak{D}_{0^+} u\right| [u]_{\Phi_x} }{[u]_{\Phi_x}}\right) d x\\
 &\leqslant[u]_{\Phi_x}^{\varphi^{-}}  \int_{\Lambda} \Phi_x\left(\tfrac{\left|\mathfrak{D}_{0^+} u\right|}{[u]_{\Phi_x}}\right) d x \\
& \leqslant[u]_{\Phi_x}^{\varphi^{-}} .
\end{aligned}
$$
Thence, we establish that $[u]_{\Phi_x}^{\varphi^{+}} \leqslant {}^0\rho_{\Phi_x}^{\psi}(u) \text{ for all } u \in \mathcal{K}_{\Phi_x} \text{ with } [u]_{\Phi_x}<1.$ Using an argument analogous to that used in the proof of \eqref{36},
, we have 
\begin{equation}\label{310}
    \Phi_x(s) \geqslant \tau^{\varphi^{+}} \Phi_x\left(\frac{s}{\tau}\right),~~ \text { for all } s>0 ~~\text { and }~~ \tau \in(0,1) \text {. }
\end{equation}
Let $u \in \mathcal{K}_{\Phi_x}$ with $[u]_{\Phi_x}<1$ and $\sigma \in\left(0,[u]_{\Phi_x}\right)$, so by $(\ref{310})$ we find
\begin{equation}\label{311}
    \int_{\Lambda} \Phi_x\left(\left|\mathfrak{D}_{0^+} u\right|\right) d x \geqslant \sigma^{\varphi^{+}} \int_{\Lambda}  \Phi_x\left(\tfrac{\left|\mathfrak{D}_{0^+} u\right|}{\sigma}\right) d x .
\end{equation}
We define $v(s)=\frac{u(s)}{\sigma}$ for all $s \in \Lambda$, we have $[v]_{\Phi_x}=\tfrac{[u]_{\Phi_x}}{\sigma}>1$.
Using \eqref{35}, we find 
\begin{equation}\label{312}
 \int_{\Lambda} \Phi_x\left(\tfrac{\left|\mathfrak{D}_{0^+} u\right|}{\sigma}\right) d x=\int_{\Lambda}  \Phi_x\left(\left|\mathfrak{D}_{0^+} v\right|\right) d x>[v]_{\Phi_x^{\varphi^{-}}}>1,
\end{equation}
Combining \eqref{311} and \eqref{312}, we deduce

$$
 \int_{\Lambda} \Phi_x\left(\left|\mathfrak{D}_{0^+} u\right|\right) d x\geqslant \sigma^{\varphi^{+}} .
$$
Letting $\sigma \nearrow[u]_{\Phi_x}$, we deduce that relation $(\ref{36})$ hold true. 
\hfill$\Box$\\

\noindent Now, for $\alpha \in (0,1)$, we assume that the function $\psi$ satisfies the condition
\begin{equation}\label{condition}
    \bigl(\psi(s)-\psi(t)\bigr)^{\alpha-1}
    < \frac{1}{\psi^{\prime}(t)}, 
    \qquad \text{for all } s \in \Lambda \text{ and } t \in [0,s].
\end{equation}

\begin{proposition}\label{P33}
Let $\alpha \in (0,1)$ and let $\Phi_x$ be a Musielak function.  
Assume that condition~\eqref{condition} holds. Then, for every $u \in L_{\Phi_x}(\Lambda)$, we have
$$
\left\{\begin{array}{l}
\left\|\mathbf{I}_{0+}^{\alpha ; \psi} u\right\|_{\Phi_x} \leqslant \left[\tfrac{\left(\psi(T)-\psi(0)\right)^{\alpha}}{\Gamma(\alpha+1)}\right]^{\frac{1}{\varphi^-}} \|u\|_{\Phi_x}^{\frac{\varphi^+}{\varphi^-}},~~ \text{ if } \|\cdot\|>1\\
\left\|\mathbf{I}_{0+}^{\alpha ; \psi} u\right\|_{\Phi_x} \leqslant \left[\tfrac{\left(\psi(T)-\psi(0)\right)^{\alpha}}{\Gamma(\alpha+1)}\right]^{\frac{1}{\varphi^+}} \|u\|_{\Phi_x}^{\frac{\varphi^-}{\varphi^+}},~~ \text{ if } \|\cdot\|<1,\\
\end{array}\right.
$$
where $\|\cdot\|>1$ means that $\|v\|_{\Phi_x}>1$ and/or $\|\mathbf{I}_{0+}^{\alpha ; \psi} v\|_{\Phi_x}>1$, The same notation is adopted for the lowercase  $(<)$.
\end{proposition}

\noindent \textbf{\textit{Proof.}} 
%Si $f$ est convexe sur un intervalle contenant 0 et $r$, alors pour tout $s \in[0,1]$
%$$f(s r) \leq s f(r)+(1-s) f(0)$$
By means of Dirichlet's formula and Jensen's inequality, we arrive at
$$\begin{aligned}
\rho_{\Phi_x}\left(\mathbf{I}_{0+}^{\alpha ; \psi} u\right)
& =\int_0^T \Phi_x\left(\left|\frac{1}{\Gamma(\alpha)} \int_0^{x} \psi^{\prime}(t)(\psi(x)-\psi(t))^{\alpha-1} u(t) d t\right|\right) dx \\
& \leqslant \int_0^T \int_0^{x} \Phi_x\left(\left|\frac{1}{\Gamma(\alpha)}  \psi^{\prime}(t)(\psi(x)-\psi(t))^{\alpha-1} u(t) d t \right|\right) dx\\
& \leqslant \frac{1}{\Gamma(\alpha)} \int_0^T \int_0^{x}  \Phi_x\left(\left|  \psi^{\prime}(t)(\psi(x)-\psi(t))^{\alpha-1} u(t)  \right|\right) d x d t  \\
& \leqslant \frac{1}{\Gamma(\alpha)} \int_0^T \int_0^{x}  \Phi_x\left(\left|u(t)  \right|\right) \psi^{\prime}(x)(\psi(x)-\psi(t))^{\alpha-1}  d x d t  \\
&= \frac{1}{\Gamma(\alpha)} \int_0^T \Phi_x\left(\left|  u(t)\right|\right) \int_{t}^{T} \psi^{\prime}(t)(\psi(x)-\psi(t))^{\alpha-1}   dx d t \\
&=\frac{\left(\psi(T)-\psi(0)\right)^{\alpha}}{\Gamma(\alpha+1)} \rho_{\Phi_x}(u).
\end{aligned}$$
Hence, if $\|\cdot\|>1$ by Proposition $\ref{Pro21}$-(i), we get 
$$\left\|\mathbf{I}_{0+}^{\alpha ; \psi} u\right\|_{\Phi_x} \leqslant \left[\frac{\left(\psi(T)-\psi(0)\right)^{\alpha}}{\Gamma(\alpha+1)}\right]^{\frac{1}{\varphi^-}}  \|u\|_{\Phi_x}^{\frac{g^+}{\varphi^-}}.$$
Similarly, if $\|\cdot\|<1$ by Proposition $\ref{Pro21}$-(ii), we have 
$$\left\|\mathbf{I}_{0+}^{\alpha ; \psi} u\right\|_{\Phi_x} \leqslant \left[\frac{\left(\psi(T)-\psi(0)\right)^{\alpha}}{\Gamma(\alpha+1)}\right]^{\frac{1}{\varphi^+}} \|u\|_{\Phi_x}^{\frac{\varphi^-}{\varphi^+}}.$$
\hfill$\Box$

\begin{proposition}\label{P34}
Let $\alpha \in (0,1)$, $\beta \in [0,1]$, and let $\Phi_x$ be a Musielak function.  
Assume that condition~\eqref{condition} is satisfied. Then,
$$
\mathbf{I}_{0+}^{\alpha ; \psi}\bigl(\mathfrak{D}_{0^+} u(t)\bigr)=u(t),
\quad \text{for all } u \in \mathcal{K}^{0}_{\Phi_x}.
$$
Moreover, the embedding $\mathcal{K}_{\Phi_x} \hookrightarrow \mathscr{C}(\Lambda)$ holds.
\end{proposition}

\noindent \textbf{\textit{Proof.}}  For any $0<t_1<t_2\leqslant T$, using \eqref{HOLDER} 
%and \eqref{23.}
, we have 
$$
	\begin{aligned}
		 \left|\mathbf{I}_{0+}^{\alpha ; \psi} u\left(t_2\right)-\mathbf{I}_{0+}^{\alpha ; \psi} u\left(t_1\right)\right| 
		& =\frac{1}{\Gamma(\alpha)}\left|\begin{array}{c}
			\int_0^{t_1} \psi^{\prime}(x)\left(\psi\left(t_2\right)-\psi(x)\right)^{\alpha-1} u(x) d x \\
			-\int_0^{t_1} \psi^{\prime}(x)\left(\psi\left(t_1\right)-\psi(x)\right)^{\alpha-1} u(x) d x \\
			+\int_{t_1}^{t_2} \psi^{\prime}(x)\left(\psi\left(t_2\right)-\psi(x)\right)^{\alpha-1} u(x) d x
		\end{array}\right| \\
        \end{aligned}
        $$

        \begin{equation}\label{35}
		\begin{aligned}
			& \leqslant \tfrac{C}{\Gamma(\alpha)}\left\|\psi^{\prime}(x)\left(\psi\left(t_2\right)-\psi(x)\right)^{\alpha-1}-\left(\psi\left(t_1\right)-\psi(x)\right)^{\alpha-1}\right\|_{\overline{\Phi}_x}\|u\|_{\Phi_x} \\
			& \quad+\left\|\psi^{\prime}(x)\left(\psi\left(t_1\right)-\psi(x)\right)^{\alpha-1}\right\|_{\overline{\Phi}_x}\|u\|_{\Phi_x} \\
			& \leqslant \tfrac{C\|u\|_{\Phi_x}}{\Gamma(\alpha)}\left[\int_0^{t_1}\overline{\Phi}_x\left( \psi^{\prime}(x)\left[\left(\psi\left(t_2\right)-\psi(x)\right)^{\alpha-1}-\left(\psi\left(t_1\right)-\psi(x)\right)^{\alpha-1}\right]\right) d x\right]^{\frac{1}{\varphi^{+}}} \\
			& \quad+\left(\int_{t_1}^{t_2}\overline{\Phi}_x\left(\psi^{\prime}(x)\left(\psi\left(t_1\right)-\psi(x)\right)^{\alpha-1}\right) d x\right)^{\frac{1}{g^{+}}} \\
			& \leqslant \tfrac{C\|u\|_{\Phi_x}\left(\overline{\Phi}_x(1)\right)^{\frac{1}{\varphi^+}}}{\Gamma(\alpha)}\left(\int_0^{t_1} \psi^{\prime}(x)\left[\left(\psi\left(t_2\right)-\psi(x)\right)^{\alpha-1}-\left(\psi\left(t_1\right)-\psi(x)\right)^{\alpha-1}\right] d x\right)^{\frac{1}{\varphi^{+}}} \\
			&\quad +\tfrac{C\|u\|_{\Phi_x}\left(\overline{\Phi}_x(1)\right)^{\frac{1}{g^+}}}{\Gamma(\alpha)} \left(\int_{t_1}^{t_2}\left(\psi^{\prime}(x)\left(\psi\left(t_2\right)-\psi(x)\right)^{\alpha-1}\right) d x\right)^{\frac{1}{g^{+}}} \\
		& \leqslant \tfrac{C\|u\|_{\Phi_x}\left(\overline{\Phi}_x(1)\right)^{\frac{1}{g^+}}}{\Gamma(\alpha)}\left(\tfrac{\left(\psi\left(t_2\right)-\psi\left(0\right)\right)^{\alpha}}{\alpha}-\tfrac{\left(\psi\left(t_2\right)-\psi(t_1)\right)^{\alpha}}{\alpha}-\tfrac{\left(\psi\left(t_1\right)-\psi(0)\right)^{\alpha}}{\alpha}\right)^{\frac{1}{\varphi^{+}}} \\
		& \quad +\tfrac{C\|u\|_{\Phi_x}\left(\overline{\Phi}_x(1)\right)^{\frac{1}{g^+}}}{\Gamma(\alpha)}\left(\tfrac{\left(\psi\left(t_2\right)-\psi\left(t_1\right)\right)^{\alpha}}{\alpha}\right)^{\frac{1}{\varphi^{+}}} \\
		& \leqslant \tfrac{4 C\left(\overline{\Phi}_x(1)\right)^{\frac{1}{\varphi^+}}\left(\psi\left(t_2\right)-\psi\left(t_1\right)\right)^{\frac{\alpha}{\varphi^+}}}{\Gamma(\alpha+1)}\|u\|_{\Phi_x} .
		\end{aligned}
	\end{equation}
Therefore
$$\left|\mathbf{I}_{0+}^{\alpha ; \psi} u\left(t_2\right)-\mathbf{I}_{0+}^{\alpha ; \psi} u\left(t_1\right)\right|\leqslant \tfrac{4 C\left(\overline{\Phi}_x(1)\right)^{\frac{1}{\varphi^+}}\left(\psi\left(t_2\right)-\psi\left(t_1\right)\right)^{\frac{\alpha}{\varphi^+}}}{\Gamma(\alpha+1)}\|u\|_{\Phi_x}. $$
From Theorem \ref{T21}, we have
	$$
	\quad \mathbf{I}_{0+}^{\alpha ; \psi}\left(\mathfrak{D}_{0^+} u(t)\right)=u(t)+C(\psi(t)-\psi(0))^{\eta-1},~~~ t \in\Lambda .
	$$
	Since $\mathfrak{D}_{0^+} u \in L_{\Phi_x}$, then by \eqref{35}, we obtain the continuity of $\mathbf{I}_{0+}^{\alpha ; \psi}\left(\mathfrak{D}_{0^+} u(t)\right)$ in $\Lambda$. As  $u(0)=0$, thus  $C=0$, which implies $$\mathbf{I}_{0+}^{\alpha ; \psi}\left(\mathfrak{D}_{0^+} u(t)\right)=u(t).$$ The result is proved.\hfill$\Box$
%******************************

\begin{proposition}\label{P35}
     Let $\alpha\in(0,1)$, $\beta\in [0,1]$ and $\Phi_x$ be a Musielak function. Assume that $(\ref{condition})$. Then for all $u\in \mathcal{K}^{0}_{\Phi_x}$, we obtain 
    \begin{equation}\label{315.}
\left\{\begin{array}{l}
\left\| u\right\|_{\Phi_x} \leqslant \left[\frac{\left(\psi(T)-\psi(0)\right)^{\alpha}}{\Gamma(\alpha+1)}\right]^{\frac{1}{\varphi^-}} \left[u\right]_{\Phi_x}^{\frac{\varphi^+}{\varphi^-}},~~ \text{ if } \|\cdot\|>1,\\
\left\| u\right\|_{\Phi_x} \leqslant \left[\frac{\left(\psi(T)-\psi(0)\right)^{\alpha}}{\Gamma(\alpha+1)}\right]^{\frac{1}{\varphi^+}} \left[ u\right]_{\Phi_x}^{\frac{\varphi^-}{\varphi^+}},~~ \text{ if } \|\cdot\|<1.\\
\end{array}\right.
\end{equation}   
  Moreover 
    \begin{equation}\label{IN}
        \|u\|_{\infty}\leqslant \frac{M\left(\overline{\Phi}_x(1)\right)^\frac{1}{\varphi+}}{\Gamma(\alpha+1)}\left(\psi(x)-\psi(0)\right)^{\frac{\alpha}{\varphi^+}} \left[ u\right]_{\Phi_x}.
        \end{equation}
\end{proposition}
\noindent \textbf{\textit{Proof.}}  Since $\mathfrak{D}_{0^+} u \in L_{\Phi_x}$ it follows from Proposition $\ref{P33}$ that 
   $$
\left\{\begin{array}{l}
\left\| \mathbf{I}_{0+}^{\alpha ; \psi} \mathfrak{D}_{0^+} u\right\|_{\Phi_x} \leqslant \left[\frac{\left(\psi(T)-\psi(0)\right)^{\alpha}}{\Gamma(\alpha+1)}\right]^{\frac{1}{\varphi^-}} \left[u\right]_{\mathcal{K}_{\Phi_x}}^{\frac{\varphi^+}{\varphi^-}},~~ \text{ if } \|\cdot\|>1,\\
\left\| \mathbf{I}_{0+}^{\alpha ; \psi} \mathfrak{D}_{0^+} u \right\|_{\Phi_x}  \leqslant \left[\frac{\left(\psi(T)-\psi(0)\right)^{\alpha}}{\Gamma(\alpha+1)}\right]^{\frac{1}{\varphi^+}} \left[ u\right]_{\mathcal{K}_{\Phi_x}}^{\frac{\varphi^-}{\varphi^+}},~~ \text{ if } \|\cdot\|<1.\\
\end{array}\right.
$$   
Using  Proposition $\ref{P34}$, we obtain the first result $(\ref{315.})$. \\
By $(\ref{HOLDER})$,
%and $(\ref{23.})$
we have for all $u\in \mathcal{K}^{0}_{\Phi_x}$ that 
$$\begin{aligned}
		\left|\mathbf{I}_{0+}^{\alpha ; \psi}\mathfrak{D}_{0^+} u\right| &=\left|\frac{1}{\Gamma(\alpha)} \int_0^x \psi^{\prime
		}(t)(\psi(x)-\psi(t))^{\alpha-1}~ \mathfrak{D}_{0^+} u(t)~dt \right|\\
		&\leqslant \frac{1}{\Gamma(\alpha)} \int_0^x \psi^{\prime
		}(t)(\psi(x)-\psi(t))^{\alpha-1}~ \left|\mathfrak{D}_{0^+} u(t)\right| ~dt\\
		&\leqslant \frac{M}{\Gamma(\alpha)} \|\psi^{\prime
		}(t)(\psi(x)-\psi(t))^{\alpha-1}\|_{\overline{\Phi}_x}\left\|\mathfrak{D}_{0^+} u\right\|_{\Phi_x}\\
		&\leqslant \frac{M}{\Gamma(\alpha)} \left(\int_0^x \overline{\Phi}_x\left(\psi^{\prime
		}(t)(\psi(x)-\psi(t))^{\alpha-1}\right) dt \right)^{\frac{1}{\varphi^+}}\left\|\mathfrak{D}_{0^+} u\right\|_{\Phi_x}\\
		&\leqslant \frac{M \left(\overline{\Phi}_x(1)\right)^{\frac{1}{\varphi^+}}}{\Gamma(\alpha+1)}\left(\psi(x)-\psi(0)\right)^{\frac{\alpha}{\varphi^+}} \left[ u\right]_{\Phi_x}.
	\end{aligned}$$
 Hence 
    $$
        \|u\|_{\infty}\leqslant \frac{M\left(\overline{\Phi}_x(1)\right)^\frac{1}{\varphi+}}{\Gamma(\alpha+1)}\left(\psi(x)-\psi(0)\right)^{\frac{\alpha}{\varphi^+}} \left[ u\right]_{\Phi_x}.
        $$\hfill$\Box$
        %***************************
         \noindent 	\begin{remark}\label{R32.} By Proposition $\ref{P35}$, we can deduce that :
\begin{enumerate}
\item[$(i)$-] In the space $\mathcal{K}^{0}_{\Phi_x}$, the norms $\|\cdot\|_{\mathcal{K}_{\Phi_x}}$ and $\left[ \cdot\right]_{\Phi_x}$ are equivalent. Therefore, we may consider $\left[ \cdot\right]_{\Phi_x}$ as a norm on this space.

\item[$(ii)$-] The space $\mathcal{K}_{\Phi_x}$ is continuously embedded into $L_{\Phi_x}(\Lambda)$.
%\item[$(iii)$-] Suppose there exists a constant $C>0$ and $t_0>0$ such that 
$$
%\Phi_x(t) \ge C\, t^{p(x)}, \quad \forall t \ge t_0.
%$$
%Then the embedding 
%$$
%\mathcal{K}_{\Phi_x} \hookrightarrow L^{p(x)}
$$ 
%is continuous.
		\end{enumerate}
	\end{remark}
    %***************************************************
 \section{Mountain Pass Solutions for $\left(\Phi_x,\psi\right)$-Fractional Musielak Equations}\label{sec4}
In this section, we investigate problem \hyperref[P]{($\mathcal{P}$)} within the new $\left(\Phi_x,\psi\right)$-fractional Musielak space $\mathcal{K}_{\Phi_x}$, together with the embedding results established in Section~\ref{sec3}. The analysis is developed through the application of the Mountain Pass Theorem, which ensures the existence of a nontrivial critical point of the corresponding energy functional under appropriate structural assumptions.

\noindent	In that context,	we say that $u \in \mathcal{K}^{0}_{\Phi_x}$ is a weak solution of problem \hyperref[P]{($\mathcal{P}$)}, if, for all $\phi \in \mathcal{K}^{0}_{\Phi_x}$, we have 

\begin{equation}\label{J}
\int_{\Lambda} a_x\left(\left|\mathfrak{D}_{0^+}  (u) \right| \right) \mathfrak{D}_{0^+} (u) \mathfrak{D}_{0^+} (\phi) dt
= \int_{\Lambda} h(t, u)\, \phi\, dt.
\end{equation}

We consider the energy functional 
$\mathfrak{J} : \mathcal{K}^{0}_{\Phi_x} \to \mathbb{R}$ 
associated with \eqref{J}, defined as follows:

$$
    \mathfrak{J}(u) = \int_{\Lambda} \Phi_x\left(\mathfrak{D}_{0^+}  u\right) \, dt 
    - \int_{\Lambda} \mathcal{H}(t, u) \, dt,
$$
where $\mathcal{H}(t, x) = \int_{0}^{x} h(t, s) \, ds.$

By a standard argument similar to that used in {\cite[Lemma 4.1]{kasmi}}, it follows that 
$\mathfrak{J}_\lambda \in \mathscr{C}^{1}\bigl(\mathcal{K}^{0}_{\Phi_x},\mathbb{R}\bigr)$, and
for all $u, \phi \in \mathcal{K}^{0}_{\Phi_x}$ we infer 
$$\left\langle \mathfrak{J}^{\prime}(u), \phi\right\rangle=\int_{\Lambda} a_x\left(\left|\mathfrak{D}_{0^+}  u \right| \right) \mathfrak{D}_{0^+}  u  \mathfrak{D}_{0^+} \phi  dt - \int_{\Lambda} h(t, u) \phi dt$$

\noindent Therefore, the weak solutions of \hyperref[P]{($\mathcal{P}$)} correspond to the critical points of the functional $\mathfrak{J}$.\\
We are now in a position to state our main existence result as follows:

\begin{theorem}\label{T41}
Let $0 < \alpha < 1$ and $0 \leq \beta \leq 1$. 
Assume that the function $h$ satisfies the conditions \hyperlink{h1}{\textbf{$($h$_1)$}} and \hyperlink{h1}{\textbf{$($h$_2)$}}. 
Then, the problem \hyperref[P]{($\mathcal{P}$)} admits at least one nontrivial weak solution 
$u \in \mathcal{K}^{0}_{\Phi_x}$.
\end{theorem}

To prove Theorem $\ref{T41}$, we make use of the Mountain Pass Theorem. 
Before proceeding, we recall and establish several auxiliary results that play a fundamental role 
in the proof of the main existence result presented in this section.
\begin{lemma}[\!\!\cite{Torres}]\label{L42}
Assume that $h$ satisfies condition \hyperlink{h2}{\textbf{(h$_2$)}}.  
Then, for every $t \in \Lambda$, the following inequalities hold:
\begin{equation}\label{43}
    \mathcal{H}(t, u) \leq \mathcal{H}\!\left(t, \frac{v}{|u|}\right)|u|^{\mu}, 
    \quad \text{for } 0 < |u| \leq 1;
\end{equation}
and
\begin{equation}\label{44}
    \mathcal{H}(t, u) \geq \mathcal{H}\!\left(t, \frac{u}{|u|}\right)|u|^{\mu}, 
    \quad \text{for } |u| \geq 1.
\end{equation}

\end{lemma}

\begin{lemma}[\!\!\cite{Torres}]\label{L43}
Let 
$\ell = \inf \{\, \mathcal{H}(t, u) \;|\; t \in \Lambda,\ |u| = 1 \,\}$.
Then, for any $s \in \mathbb{R} \setminus \{0\}$ and any $u \in \mathcal{K}^{0}_{\Phi_x}$, we have
\begin{equation}\label{45}
    \int_{\Lambda} \mathcal{H}(t, s u(t))\, dt 
    \geq \ell\, |s|^{\mu} \int_{\Lambda} |u(t)|^{\mu}\, dt 
    - T\ell.
\end{equation}
\end{lemma}

The following result, commonly known as the Palais-Smale compactness condition (PS), provides the compactness criterion essential for the application of the Mountain Pass Theorem.

\begin{lemma}
Let $\Phi_x$ be a Musielak function satisfying $(\ref{delta})$, and assume that the function $h$ fulfills the conditions \hyperlink{h1}{\textbf{$($h$_1)$}} and \hyperlink{h2}{\textbf{$($h$_2)$}}. Then, the functional $\mathfrak{J}$ verifies the Palais-Smale condition.
\end{lemma}

\noindent \textbf{\textit{Proof.}} Let $\left\{u_k\right\}$ be a (PS)-sequence of $\mathfrak{J}$ on $\mathcal{K}^{0}_{\Phi_x}$, which could be expressed mathematically as follows
\begin{equation}\label{45}
    \left|\mathfrak{J}\left(u_k\right)\right| \leq M\quad \text{and}\quad \lim _{k \rightarrow \infty} \mathfrak{J}^{\prime}\left(u_k\right)=0 .
\end{equation}
We begin by showing that the sequence $\{u_k\}$ is bounded, recall that
$$ \mathfrak{J}(u_k) = \int_{\Lambda} \Phi_x\left(\mathfrak{D}_{0^+}  u_n\right) \, dt 
    - \int_{\Lambda} \mathcal{H}(t, u_k) \, dt,$$
and 
$$
 \left\langle \mathfrak{J}^{\prime}(u_k), u
 _k\right\rangle=\int_{\Lambda} a_x\left(\left|\mathfrak{D}_{0^+}  u_k \right| \right) \mathfrak{D}_{0^+}  u_k  \mathfrak{D}_{0^+} u
 _k  dt - \int_{\Lambda} h(t, u
 _k) u
 _k dt.
$$
Then by $(\ref{45})$, we get
\begin{equation}
    \begin{aligned}\left|\mathfrak{J}\left(u_k\right)-\tfrac{1}{\mu} \left\langle \mathfrak{J}^{\prime}\left(u_k\right), u_k\right\rangle \right| & \leq\left|\mathfrak{J}\left(u_k\right)\right|+\left|\tfrac{1}{\mu} \mathfrak{J}^{\prime}\left(v_k\right)\right|\left|v_k\right| \\ & \leq C\left(1+\left\|u_k\right\|_{\mathcal{K}^{0}_{\Phi_x}}\right) .\end{aligned}
\end{equation}
However, with Proposition $\ref{P32}$ and $(\ref{S})$, we get

\begin{equation}\label{S6}
    \begin{aligned} 
\mathfrak{J}( u_k) & = \int_{\Lambda} \Phi_x\left(\mathfrak{D}_{0^+}  u_k\right) \, dt 
    - \int_{\Lambda} \mathcal{H}(t, u_k) \, dt \\ 
& \geqslant \left\|u_k\right\|_{\mathcal{K}^{0}_{\Phi_x}}^{\varphi^\pm}- \int_{\Lambda} \mathcal{H}(t, u_k(t))  d t.
\end{aligned}
\end{equation}
In addition
\begin{equation}\label{S7}
    \begin{aligned} 
   \left\langle  \mathfrak{J}^{\prime}(u_k), u_k\right\rangle & =\int_{\Lambda} a_x\left(\left|\mathfrak{D}_{0^+}  u_k \right| \right) \mathfrak{D}_{0^+}  u_k  \mathfrak{D}_{0^+} u_k  dt - \int_{\Lambda} h(t, u_k) u_k dt \\ & \leqslant k \int_{\Lambda} \Phi_x\left(\mathfrak{D}_{0^+}  u_k\right) dt - \int_{\Lambda} h(t, u_k(t)) u_k(t) d t\\
    &\leqslant k \left\|u_k\right\|_{\mathcal{K}^{0}_{\Phi_x}}^{\varphi^\pm} - \int_{\Lambda} h(t, u_k(t)) u_k(t) d t.
    \end{aligned}
    \end{equation}
By $(\ref{S6})$, $(\ref{S7})$ and \hyperlink{h2}{\textbf{$($h$_2)$}}, we have
\begin{equation}
    \begin{aligned}  \mathfrak{J}\left(u_k\right)-\tfrac{1}{\mu} \left\langle \mathfrak{J}^{\prime}\left(u_k\right), u_k\right\rangle   \geqslant & \left(1-\tfrac{k}{\mu}\right)\left\|u_k\right\|_{\mathcal{K}^{0}_{\Phi_x}}^{\varphi^\pm}-\int_{\Lambda} \mathcal{H}\left(t, u_k(t)\right) d t+\tfrac{1}{\mu} \int_{\Lambda} h\left(t, u_k(t)\right) u_k(t) d t \\
    \geq & \left(1-\tfrac{k}{\mu}\right)\left\|u_k\right\|_{\mathcal{K}^{0}_{\Phi_x}}^{\varphi^\pm} .
    \end{aligned}
    \end{equation}
Since $\mu > k$, it follows that the sequence $\{v_k\}$ is bounded in $\mathcal{K}^{0}_{\Phi_x}$. 
As $\mathcal{K}^{0}_{\Phi_x}$ is a reflexive space, there exists a function $u \in \mathcal{K}^{0}_{\Phi_x}$ 
and a subsequence, still denoted by $\{u_k\}$ for simplicity, such that
$$
u_k \rightharpoonup u \quad \text{in } \mathcal{K}^{0}_{\Phi_x} \text{ as } k \to \infty.
$$

\begin{equation}\label{410}
    \begin{aligned} 
\left\langle\mathfrak{J}^{\prime}\left(u_k\right)-\mathfrak{J}^{\prime}(u), u_k-v\right\rangle  = & \left\langle\mathfrak{J}^{\prime}\left(u_k\right), u_k-u\right\rangle-\left\langle \mathfrak{J}^{\prime}(u), u_k-u\right\rangle \\ \leq & \left\|\mathfrak{J}^{\prime}\left(v_k\right)\right\|\left\|u_k-u\right\|_{\mathcal{K}^{0}_{\Phi_x}}-\left\langle\mathfrak{J}^{\prime}(u), u_k-u\right\rangle .
\end{aligned}
\end{equation}
By $(\ref{410})$, we find 
$$\left\langle\mathfrak{J}^{\prime}\left(u_k\right)-\mathfrak{J}^{\prime}(u), u_k-u\right\rangle \rightarrow 0, \text{ as } k \rightarrow \infty.
$$
From Propositions $\ref{P34}$ and $\ref{P35}$, we get that $u_k$ is bounded in $\mathscr{C}(\Lambda)$, additionally, we can reasonably suppose that
$$\lim _{k \rightarrow \infty}\left\|u_k-u\right\|_{\infty}=0.$$
Hence, we deduce
$$
\int_{\Lambda}\left[h\left(t, u_k(t)\right)-h(t, u(t))\right]\left(u_k(t)-u(t)\right) d t \underset{k \rightarrow \infty}{\longrightarrow}0.
$$
Moreover, a straightforward calculation shows that
$$
\left\langle \mathfrak{J}^{\prime}\left(u_k\right)-\mathfrak{J}^{\prime}(u), u_k-u\right\rangle\geqslant\left\|u_k-u\right\|_{\mathcal{K}^{0}_{\Phi_x}}^{\varphi^\pm}-\int_{\Lambda}\left(h\left(t, u_k(t)\right)-h(t, u(t))\right)\left(u_k(t)-u(t)\right) d t .
$$

%1. For almost every $x \in \Lambda$, the function $s \mapsto a_x(|s|)s$ is monotone (strictly monotone if necessary), i.e.

%$$
%\left(a_x(|\xi|)\xi - a_x(|\eta|)\eta\right) \cdot (\xi - \eta) \geq 0 \quad \text{for all } \xi, \eta \in \mathbb{R}^N.
%$$

%2. A coercivity/modularity condition links the operator to the Young function $\boldsymbol{\Phi}_x$: there exists a constant $c > 0$ such that, for almost every $x$ and all $\xi, \eta$,
%$$
%\left(a_x(|\xi|)\xi - a_x(|\eta|)\eta\right) \cdot (\xi - %\eta) \geq c\, \Phi_x(|\xi - \eta|).
%$$

%\begin{aligned}
%&\begin{aligned}
%& \int_{\Lambda} \left[a_x\left(\left|\mathfrak{D}_{0^{+}} v_k\right|\right)\mathfrak{D}_{0^{+}} v_k 
%- a_x\left(\left|\mathfrak{D}_{0^{+}} v\right|\right)\mathfrak{D}_{0^{+}} v \right] 
%\mathfrak{D}_{0^{+}}\left(v_k - v\right) \, dx \\
%& \quad \geq c \int_{\Lambda} \Phi_x\left(\left|\mathfrak{D}_{0^{+}} v_k - \mathfrak{D}_{0^{+}} v\right|\right) \, dx 
%= c\, \rho\left(v_k - v\right)
%\end{aligned} \\
%&\text{By combining this with the nonlinear term, we obtain:} \\
%&\left\langle \mathfrak{J}^{\prime}\left(v_k\right) - \mathfrak{J}^{\prime}(v), v_k - v \right\rangle 
%\geq c \int_{\Lambda} \Phi_x\left(\left|\mathfrak{D}%_{0^{+}} v_k - \mathfrak{D}_{0^{+}} v\right|\right) \, dx 
%- \int_{\Lambda} \left[h\left(t, v_k\right) - h(t, v)\right]\left(v_k - v\right) \, dt
%\end{aligned}

Hence, 
$$
\|u_k - u\|_{\mathcal{K}^{0}_{\Phi_x}} \longrightarrow 0 \quad \text{as } k \to \infty,
$$
which means that the sequence $\{u_k\}$ converges strongly to $u$ in $\mathcal{K}^{0}_{\Phi_x}$.
\hfill$\Box$\\

\noindent We are now ready to establish the proof of the main result of this section.

\noindent \textbf{\textit{Proof of Theorem}} $\ref{T41}$. Evidently, $\mathfrak{J}(0)=0$. It thus remains to establish that $\mathfrak{J}$ meets the geometric criteria stipulated by the mountain pass theorem.\\

From $(\ref{IN})$, we obtain
$$
\max_{t \in \Lambda} |u(t)| \le R \|u\|_{\mathcal{K}^{0}_{\Phi_x}}, \quad \forall u \in \mathcal{K}^{0}_{\Phi_x},
$$
where
$$
R = \frac{M\left(\overline{\Phi}_x(1)\right)^\frac{1}{\varphi+}}{\Gamma(\alpha+1)}\left(\psi(x)-\psi(0)\right)^{\frac{\alpha}{\varphi^+}}.
$$
Next, let $C_1 = \frac{1}{R}$. Then, by the above inequality together with $(\ref{43})$, we deduce that if $\|u\|_{\mathcal{K}^{0}_{\Phi_x}} \le C_1$, 
$$
\begin{aligned}
\int_{\Lambda} \mathcal{H}(t, u(t)) \, dt 
& \le \int_{\Lambda} \mathcal{H}\Big(t, \frac{u(t)}{|v(t)|}\Big) |u(t)|^\mu \, dt \\
& \le R^\mu T \tilde{C} \, \|u\|_{\mathcal{K}^{0}_{\Phi_x}}^\mu.
\end{aligned}
$$

Then
$$\begin{aligned} \mathfrak{J}( u) & = \int_{\Lambda} \Phi_x\left(\mathfrak{D}_{0^+}  u\right) \, dt 
    - \int_{\Lambda} \mathcal{H}(t, u) dt\,
\\ & \geq \|u\|_{\mathcal{K}^{0}_{\Phi_x}}^{\varphi^\pm}-R^{\mu}T \tilde{C}\|u\|_{\mathcal{K}^{0}_{\Phi_x}}^{\mu}, \quad \text { if }\|u\|_{\mathcal{K}^{0}_{\Phi_x}} \leq C_1,\end{aligned}$$

Therefore
$$\mathfrak{J}(v)\geq C_1^{\varphi^{\pm}}-R^{\mu}TC_1^{\mu} \tilde{C}, \text{ if } |v\|_{\mathcal{K}^{0}_{\Phi_x}} = C_1. $$

Let us consider $L<\min \left\{C_1,\left(\tfrac{1}{ R^\mu T \tilde{C}}\right)^{\frac{1}{\mu-\varphi^{\pm}}}\right\}$ and $\theta=L^{\varphi^\pm}-L^\mu T \Tilde{C} R^\mu$, then
$$
\mathfrak{J}(u) \geq \theta ~~\text { with }\|u\|_{\mathcal{K}^{0}_{\Phi_x}}=L.
$$
Therefore, $\mathfrak{J}$ satisfies the first geometric condition required by the mountain pass theorem.\\

Let $s \in \mathbb{R} \setminus \{0\}$ and $u \in \mathcal{K}^{0}_{\Phi_x}$. From Lemma~\ref{L43}, we drive 
$$
\mathfrak{J}(s u) \le s \, \|u\|_{\mathcal{K}^{0}_{\Phi_x}}^{\varphi^\pm} - \ell |s|^\mu \int_{\Lambda} |u(t)|^\mu \, dt + T \ell.
$$  
Since $\mu > k$, letting $s \to +\infty$ gives $\mathfrak{J}(s u) \to -\infty$.  
Hence, the second geometric condition of the mountain pass theorem is satisfied by taking $e = s v$ with $s$ sufficiently large, so that $\mathfrak{J}(e) \le 0$.

\noindent Therefore, the functional $\mathfrak{J}$ satisfies the mountain pass condition. Consequently, $\mathfrak{J}$ admits a nontrivial critical point, which corresponds to a nontrivial weak solution of problem \hyperref[P]{($\mathcal{P}$)}.\\

\noindent\textbf{{\large Declarations}} :\\

\noindent\textbf{Ethical Approval :} Not applicable.\\

\noindent\textbf{Competing interests :} The authors declare that there is no conflict of interest.\\

\noindent\textbf{Authors' contributions :} The authors contributed equally to this work.\\

\noindent\textbf{Funding :} Not applicable.\\

\noindent\textbf{Availability of data and materials :} Not applicable.\\

    % % % % % % % % % % % % % % % % % % % % % % % % % % % % % % % % % % % % % % % % % % % % % % % % % % % % % % % % % % % % % % % % % % % % % % % % % % % % % % % % % % % % % % % % % % % % % % % %

\end{document}